\documentclass{amsart}
\usepackage{amssymb}

\newtheorem{thm}{Theorem}[section]
\newtheorem{prop}[thm]{Proposition}
\newtheorem{cor}[thm]{Corollary}
\newtheorem{lem}[thm]{Lemma}

\theoremstyle{definition}
\newtheorem{defn}[thm]{Definition}

\theoremstyle{remark}

\newtheorem{rem}[thm]{Remark}
\newtheorem{example}[thm]{Example}

\newtheorem*{ack}{Acknowledgment}

\newcommand{\GS }{\Gamma (G, S)}
\newcommand{\HT }{\Gamma (H, T)}

\usepackage[bookmarks=true, pdfauthor={YANG Wenyuan}]{hyperref}

\begin{document}

\title{limit sets of relatively hyperbolic groups}

\author{Wen-yuan Yang}
\address{College of Mathematics and Econometrics, Hunan
University, Changsha, Hunan 410082 People's Republic of China}
\curraddr{U.F.R. de Mathematiques, Universite de Lille 1, 59655
Villeneuve D'Ascq Cedex, France}
\email{wyang@math.univ-lille1.fr}
\thanks{The author is supported by the China-funded Postgraduates
Studying Aboard Program for Building Top University. This research
was supported  by National Natural Science Foundational of China
(No. 11071059)}



\date{}

\dedicatory{}


\begin{abstract}
In this paper, we prove a limit set intersection theorem in
relatively hyperbolic groups. Our approach is based on a study of
dynamical quasiconvexity of relatively quasiconvex subgroups. Using
dynamical quasiconvexity, many well-known results on limit sets of
geometrically finite Kleinian groups are derived in general
convergence groups. We also establish dynamical quasiconvexity of
undistorted subgroups in finitely generated groups with nontrivial
Floyd boundary.
\end{abstract}

\maketitle

\section{Introduction}
The purpose of this paper is to study applications of dynamical
quasiconvexity to limit sets of relatively quasiconvex subgroups.
The notion of dynamical quasiconvexity is introduced by Bowditch in
\cite{Bow2} and used to characterize a geometric notion of
quasiconvexity in word hyperbolic groups. In relatively hyperbolic
groups, Gerasimov-Potyagailo \cite{GePo1} recently showed that
dynamically quasiconvex subgroups are exactly the class of
relatively quasiconvex subgroups.

In this paper, we shall show that even in general convergence
groups, dynamical quasiconvexity efficiently captures algebraic and
geometric properties of subgroups as well. However, the particular
interest we have in mind lies in relatively hyperbolic groups, and
finitely generated groups with nontrivial Floyd boundary, which is
conjectured to be relatively hyperbolic \cite{OOS}.

Let $G$ be a finitely generated group, admitting a
\textit{convergence group action} on a compact metric space $M$.
Then the \textit{limit set} $\Lambda(H)$ of a subgroup $H \subset G$
is the set of accumulation points of $H$-orbits in $M$. See Section
2 for their precise definitions. Following Anderson \cite{Aj}, a
limit set intersection theorem for convergence groups describes the
limit set $\Lambda(H \cap J)$ in terms of $\Lambda(H)$ and
$\Lambda(J)$, where $H$, $J$ are subgroups of a convergence group
$G$. Ideally, we expect such a theorem has the following form
\begin{center}
$\Lambda(H) \cap \Lambda(J) = \Lambda(H \cap J) \cup E$\\
\end{center}
where $E$ is an exceptional set consisting of specific parabolic
points of $\Lambda(H)$ and $\Lambda(J)$.

Such a limit set intersection theorem has been investigated in
several different classes of groups. In 1992, Susskind-Swarup
\cite{SuSw} showed that the above decomposition of limit sets holds
for a pair of geometrically finite Kleinian subgroups. In
\cite{Aj3}, \cite{Aj2} and \cite{Aj}, using techniques specific to 3
manifolds, Anderson carried out a systematic study of the
intersection of two finitely generated subgroups of 3 dimensional
Kleinian groups and proved that the limit set intersection theorem
holds in this context.

In 1987, Gromov \cite{Gro} introduced relatively hyperbolic groups
as a generalization of many naturally occurred groups, for example,
word hyperbolic groups and geometrically finite Kleinian groups and
many others. In word hyperbolic groups, the limit set intersection
theorem is explained in Gromov \cite[Page 164]{Gro2}, where the
exceptional set $E$ is empty. Our main result is to generalize these
limit set intersection theorems in relatively hyperbolic groups as
follow.

\begin{thm} \label{thm:intersection}
Let $H, J$ be two relatively quasiconvex subgroups of a relatively
hyperbolic group $G$. Then
\begin{center}
$\Lambda(H) \cap \Lambda(J) = \Lambda(H \cap J) \bigsqcup E$\\
\end{center}
where the exceptional set $E$ consists of parabolic fixed points of
$\Lambda(H)$ and $ \Lambda(J)$, whose stabilizer subgroups in $H$
and $J$ have finite intersection. Equivalently, the set $E$ consists
of the limit points isolated in $\Lambda(H) \cap \Lambda(J)$.
\end{thm}

\begin{rem}
Two special cases of Theorem \ref{thm:intersection} were known in
Dahmani \cite{Dah}. The first one is under the assumption that
maximal parabolic groups are abelian. The other one is proved for a
pair of fully quasiconvex subgroups, where the exceptional set $E$
is empty. By a result of Hruska \cite{Hru}, fully quasiconvex
subgroups are relatively quasiconvex .
\end{rem}

One corollary to Theorem \ref{thm:intersection} is the following
well-known result, which is usually proved via geometrical
methods by Hruska \cite{Hru} and, independently, Martinez-Pedroza \cite{MarPed}.

\begin{cor}
Let $H, J$ be two relatively quasiconvex subgroups of a relatively
hyperbolic group $G$. Then $H \cap J$ is relatively quasiconvex.
\end{cor}

The proof of Theorem \ref{thm:intersection} replies crucially on
dynamical quasiconvexity of relatively quasiconvex subgroups. Even
in general convergence groups, dynamical quasiconvex subgroups share
many nice properties with relatively quasiconvex subgroups. See
Section 2 for a few of them.

The following result extends a property of quasiconvex subgroups in
word hyperbolic groups proved in Mihalik-Towle \cite{MM} to
dynamically quasiconvex subgroups in general convergence groups.
\begin{thm} \label{thm:conjuagte}
Let $H$ be dynamically quasiconvex in a convergence group $G$ such
that $|\Lambda(H)| \geq 2$. Then for any $g \in G \setminus H$,
$gHg^{-1} \subseteq H$ implies that $gHg^{-1} = H$.
\end{thm}

In the final section, we explore dynamical quasiconvexity of
subgroups in a special class of convergence groups, i.e. finitely
generated groups with nontrivial Floyd boundary. In \cite{Floyd},
Floyd boundary was introduced by Floyd to compactify Cayley graphs
of finitely generated groups. Later, Karlsson \cite{Ka} proved that
the left multiplication of group elements extends to a convergence
group action on Floyd boundary.

Our last result establishes the dynamical quasiconvexity of
undistorted subgroups with respect to the convergence action on
Floyd boundary.
\begin{thm} \label{thm:dynconvex}
If $H$ is an undistorted subgroup of a finitely generated group $G$
with nontrivial Floyd boundary, then $H$ is dynamical quasiconvex.
\end{thm}
\begin{rem}
The class of groups with nontrivial Floyd boundary includes
non-elementary finitely-generated relatively hyperbolic groups
\cite{Ge2}. In particular, infinite-ended groups have nontrivial
Floyd boundary.
\end{rem}

By the Floyd map theorem in \cite{Ge2}, there exists an equivariant
map from the Floyd boundary of a relatively hyperbolic group to its
Bowditch boundary. Then it is easily seen that the dynamical
quasiconvexity of a subgroup is kept under an equivalent quotient.
See a proof in \cite[Lemma 4.5]{YANG3}, for example. This together
with Theorem \ref{thm:dynconvex} gives the following result, which
is first proved in \cite{Hru}.

\begin{cor}
Let $H$ be an undistorted subgroup of a relatively hyperbolic group
$G$. Then $H$ is relatively quasiconvex.
\end{cor}

The paper is organized as follows. In Section 2, we define dynamical
quasiconvex subgroups in general convergence groups. Then we deduce
several consequences of dynamically quasiconvex subgroups in general
convergence groups. In Section 3, we study the intersection of
conical limit points and bounded parabolic points of dynamically
quasiconvex subgroups, and then conclude with the proof of Theorem
\ref{thm:intersection}. In Section 4, it is shown that a
nonparabolic dynamically quasiconvex subgroup cannot contain a
proper conjugate of itself. In the final section, we briefly
introduce the Floyd boundary for a finitely generated groups and
give a proof of Theorem \ref{thm:dynconvex}. Moreover, an example is
given of dynamically quasiconvex subgroups which are not
geometrically finite.

\begin{ack}
The author would like to sincerely thank Prof. Leonid Potyagailo for
many helpful comments and inspired discussions during the course of
this work. The author also thanks Chris Hruska for pointing out
several inaccuracies of references. The author is grateful for many
valuable comments of the referee, which substantially improves the
presentation of the paper.
\end{ack}

\begin{center} \section{Preliminary results}\end{center}

Throughout the paper, we consider a finitely generated group $G$,
and a compact metrizable space $M$ containing at least three points.

A \textit{convergence group action} is an action of a group $G$ on
$M$ such that the induced action of $G$ on the space $\Theta M$ of
distinct unordered triples of points of $M$ is properly
discontinuous.

Suppose $G$ has a convergence group action on $M$. Then $M$ is
partitioned into a limit set $\Lambda(G)$ and discontinuous domain
$M \setminus \Lambda(G)$. The \textit{limit set} $\Lambda(H)$ of a
subgroup $H \subset G$ is the set of limit points, where a
\textit{limit point} is an accumulation point of some $H$-orbit in
$M$. It is well-known that if $|\Lambda(H)| \ge 2$, the limit set
$\Lambda(H)$ is also characterized as the minimal $H$-invariant
closed subset in $M$ of cardinality at least two.

An element $g \in G$ is \textit{elliptic} if it has finite order. An
element $g \in G$ is \textit{parabolic} if it has infinite order and
fixes exactly one point of $M$. An element $g \in G$ is
\textit{loxodromic} if it has infinite order and fixes exactly two
points of $M$. An infinite subgroup $P \subset G$ is a
\textit{parabolic subgroup} if it contains no loxodromic element. A
parabolic subgroup $P$ has a unique fixed point in $M$. This point
is called a \textit{parabolic point}. The stabilizer of a parabolic
point is always a maximal parabolic group. A parabolic point $p$
with stabilizer $G_p$ := $Stab_G (p)$ is \textit{bounded} if $G_p$
acts cocompactly on $M \setminus \{p\}$. A point $z \in M$ is a
\textit{conical limit} point if there exists a sequence $\{g_i\}$ in
$G$ and distinct points $a, b \in M$ such that $g_i (z) \to a$ ,
while for all $q \in M \setminus \{z\}$ we have $g_i (q) \to b$.

Before discussing relative hyperbolicity, we recall the following
well-known result on general convergence groups.

\begin{lem} \cite[Theorem 3.A]{Tukia} \label{lem:conipara}
In a convergence group, a conical limit point can not be parabolic.
\end{lem}

In the literature, various definitions of relative hyperbolicity
were proposed, see Farb \cite{Farb}, Bowditch \cite{Bow1}, Osin
\cite{Osin}, Drutu-Sapir \cite{DruSapir} and so on. These different
definitions are now proven to be equivalent for countable groups,
see Hruska \cite{Hru} for a complete account, and provide convenient
and complement viewpoints to the study of this class of groups. For
the sake of the purpose of this paper, we use the following
dynamical formulation of relatively hyperbolic groups.

\begin{defn}\label{def:RH1}
A convergence group action of $G$ on $M$ is \textit{geometrically
finite} if every limit point of $G$ is either conical or bounded
parabolic. Let $\mathbb{P}$ be a set of representatives of the
conjugacy classes of maximal parabolic subgroups. Then we say the
pair (G, $\mathbb{P}$) is \textit{relatively hyperbolic}. When
$\mathbb P$ is clear in the context, we just say $G$ is relatively
hyperbolic.
\end{defn}

\begin{rem} \label{rem:undistorted}
By the work of Drutu-Sapir \cite{DruSapir}, Osin \cite{Osin} and
Gerasimov \cite{Ge1}, maximal parabolic subgroups are quasiconvex
and finitely generated. So in the definition of relative
hyperbolicity, we do not need impose the ``finitely generated''
condition on maximal parabolic subgroups, as usually do in Bowditch
\cite{Bow1}.
\end{rem}

The following notion of dynamical convexity is introduced by
Bowditch \cite{Bow2} and proven to be equivalent to the geometric
quasiconvexity in word hyperbolic groups.

\begin{defn} \label{def:DQC}
A subgroup $H$ of a convergence group $G$ is \textit{dynamically
quasiconvex} if the following set
\begin{center}
$\{gH \in G/H : g \Lambda(H) \cap K \neq \emptyset$\;and\;$g
\Lambda(H) \cap L \neq \emptyset\}$
\end{center}
is finite, whenever $K$ and $L$ are disjoint closed subsets of $M$.
\end{defn}

Recently, Gerasimov-Potyagailo \cite{GePo1} proved that dynamical
quasiconvexity coincides with relatively quasiconvexity in
relatively hyperbolic groups, which gives a positive answer to a
question of Osin in his book \cite{Osin}. We refer the reader to
\cite{Osin} for the definition of relative quasiconvexity.

\begin{thm} \label{thm:DQCRQC} \cite{GePo1}
Suppose $G$ is relatively hyperbolic. Every subgroup $H$ of $G$ is
dynamically quasiconvex if and only if it is relatively quasiconvex.
\end{thm}

Let us first draw some consequences of dynamical quasiconvexity,
without assuming $G$ is relatively hyperbolic.

\begin{lem}\label{lem:samelimitset}
Let $H$ be dynamically quasiconvex in a convergence group $G$ such
that $|\Lambda(H)| \geq 2$. Then for any subgroup $H \subset J
\subset G$ satisfying $\Lambda(H)=\Lambda(J)$, we have that $H$ is
of finite index in $J$. In particular, $J$ is dynamically
quasiconvex.
\end{lem}
\begin{proof}
Since $|\Lambda(H)| \geq 2$, we can pick distinct points $x$ and $y$
from $\Lambda(H)$. Since $\Lambda(H)=\Lambda(J)$, we have that each
coset of $H$ in $J$ belongs to the following set
\begin{center}
$\{gH \in G/H : g \Lambda(H) \cap \{x\} \neq \emptyset$ and $g
\Lambda(H) \cap \{y\} \neq \emptyset\}$.
\end{center}
By the dynamical quasiconvexity of $H$, the above set is finite.
Thus, $H$ is of finite index in $J$.

In order to prove the dynamical quasiconvexity of $J$, it suffices
to show that the following set
\begin{center}
$\Omega := \{gJ \in G/J : g \Lambda(J) \cap K \neq \emptyset$ and $g
\Lambda(J) \cap L \neq \emptyset\}$
\end{center}
is finite, for any given disjoint closed subsets $K, L \subset M$. On
the other hand, since $H$ is dynamically quasiconvex, the following
set
\begin{center}
$\{gH \in G/H : g \Lambda(H) \cap K \neq \emptyset$ and $g
\Lambda(H) \cap L \neq \emptyset\}$.
\end{center}
is finite. Combining with the fact $H$ is of finite index in $J$, we
have that the above set $\Omega$ is finite. Therefore, $J$ is
dynamically quasiconvex.
\end{proof}

\begin{cor}
Let $H$, $J$ be dynamically quasiconvex in a convergence group $G$
such that $\Lambda(H) = \Lambda(J)$. If $|\Lambda(H)| \geq 2$, then
$H$ and $J$ are commensurable.
\end{cor}
\begin{proof}
Let $L$ be the stabilizer in $G$ of the limit set $\Lambda(H)$.
Using Lemma \ref{lem:samelimitset}, we have that $H, J$ are both of
finite index in $L$. It thus follows that $H\cap J$ is of finite
index in both $H$ and $J$.
\end{proof}

Recall that the \textit{commensurator} of $H$ in a convergence group
$G$ is defined as the subgroup of $G$, consisting of all $g \in G$
such that $H \cap gHg^{-1}$ has finite index in both $H$ and
$gHg^{-1}$.

\begin{cor} \label{cor:commensurator}
Let $H$ be dynamically quasiconvex in a convergence group $G$ such
that $|\Lambda(H)| \geq 2$. Then $H$ is of finite index in its
commensurator.  In particular, $H$ is of finite index in its
normalizer.
\end{cor}
\begin{proof}
Let $L$ be the commensurator of $H$ in $G$. Then $H \subset L$. It
is obvious that $\Lambda(H) \subset \Lambda(L)$.

It is well-known that the limit set of a subgroup is same as the one
of its finite extension. So for each $g \in L$, we have $\Lambda(H)
= \Lambda(gHg^{-1}) = g\Lambda(H)$, i.e. $L$ leaves invariant the
limit set $\Lambda(H)$. Since the limit set $\Lambda(L)$ is the
minimal $L$-invariant closed subset in $M$ of cardinality at least
two, we have $\Lambda(L)=\Lambda(H)$. The conclusion now follows
from Lemma \ref{lem:samelimitset}.
\end{proof}
\begin{rem}
In relatively hyperbolic groups, Corollary \ref{cor:commensurator}
has been proven using different methods in Hruska-Wise \cite{HrWi}.
We remark the hypothesis $|\Lambda(H)| \ge 2$ is necessary for the
above lemma and corollaries, as it is easy to get counterexamples
when $H$ is taken as a parabolic subgroup.
\end{rem}

\begin{center} \section{Intersections of limit sets} \end{center}
In this section, we study the intersection of limit sets of
dynamically quasiconvex subgroups. The intersection of conical limit
points is firstly examined.

\begin{prop} \label{prop:conical}
Let $H$ be dynamically quasiconvex in a convergence group $G$.
Suppose $J < G$ is infinite and let $z \in \Lambda(H) \cap
\Lambda(J)$ be a conical limit point of $J$. Then $z \in \Lambda(H
\cap J)$ is a conical limit point of $H\cap J$.
\end{prop}

\begin{proof}
Let $z \in \Lambda(H) \cap \Lambda(J)$ be a conical limit point of
$J$. Then there exists a sequence $\{j_n\}$ in $J$ and distinct
points $a, b \in \Lambda(J)$ such that $j_n (z) \rightarrow a$,
while $j_n (q)\rightarrow b$ for all $q \in \Lambda(J) \setminus
\{z\}$. By the convergence property of $\{j_n\}$, we also have that
$j_n (q) \rightarrow b$ for all $q \in M \setminus \{z\}$. In
particular, we can choose $q$ to be a limit point in $\Lambda(H)
\setminus \{z\}$. Here, we use the fact $|\Lambda(H)| \geq 2$, which
follows from Lemma \ref{lem:conipara}.

Take closed neighborhoods $U$ and $V$ of $a$ and $b$ respectively,
such that $U \cap V = \emptyset$. After passage to a subsequence of
$\{j_n\}$, we can assume $j_n(z) \in U$ and $j_n(q) \in V$ for all
$n$. This implies that $j_n H$ belongs to the following set for all
$n$,
\begin{center}
$\{gH \in G/H : g \Lambda(H) \cap U \neq \emptyset$ and $g
\Lambda(H) \cap V \neq \emptyset\}.$
\end{center}
By the dynamical quasiconvexity of $H$ in $G$, the above set is
finite. Thus, $\{j_n H\}$ is a finite set of cosets. By taking
further a subsequence of $\{j_n\}$, we suppose $j_n H = j_1 H$ for
all $n$. We can write $j_n = j_1 h_n$ for each $n$, where $h_n \in
H$. Then $j_1^{-1}j_n = h_n$ implies that $H\cap J$ is nontrivial
and infinite.

It suffices to prove that $z$ is a conical limit point of $H \cap
J$. By the convergence property of $\{j_n\}$, it follows that
$h_n(z) = j_1^{-1}j_n(z) \rightarrow j_1 ^{-1}(a)$ and $h_n(q) = j_1
^{-1}j_n(q) \rightarrow j_1 ^{-1}(b)$ for all $q \in M \setminus
\{z\}$. Thus, $z$ is a conical limit point of $H \cap J$.
\end{proof}

\begin{rem}
A similar statement of Proposition \ref{prop:conical} in relatively
hyperbolic groups appears in the proof of Proposition 3.1.10 in
\cite{Dah}.
\end{rem}

We now study how bounded parabolic points intersect. Compared to
that of conical points, the intersection of bounded parabolic points
raises some complicated behavior.

\begin{prop} \label{prop:boundpara}
Let $H$, $J$ be infinite subgroups of a convergence group $G$. If $z
\in \Lambda(H) \cap \Lambda(J)$ is a bounded parabolic point of $H$
and $J$, then $z$ is either a bounded parabolic point of $H\cap J$,
or an isolated point in $\Lambda(H) \cap \Lambda(J)$ and does not
lie in $\Lambda(H \cap J)$.
\end{prop}

\begin{proof}
Since $z$ is a bounded parabolic point of both $H$ and $J$, there
are compact subsets $K \subset \Lambda(H)\setminus z$ and $L \subset
\Lambda(J)\setminus \{z\}$, such that $H_z K=\Lambda(H)\setminus
\{z\}$ and $J_z L=\Lambda(J)\setminus \{z\}$. Here, $H_z$ and $J_z$
are stabilizers in $H$ and $J$ of $z$, respectively. Let $P = H_z
\cap J_z$.

We claim that there exists a compact subset $C \subset M \setminus
\{z\}$ such that $\Lambda(H \cap J)\setminus z \subset P C$.

Note first that $\Lambda(H \cap J)\setminus \{z\} \subset
(\Lambda(H) \cap \Lambda(J)) \setminus z = H_z K\cap J_z L$.
Therefore, it suffices to show that there exists a compact subset $C
\subset M\setminus \{z\}$ such that $H_z K\cap J_z L \subset P C$.
Since $G$ is countable, we define the following set
\begin{equation}\label{countset}
\mathcal{A} = \{h_n K\cap j_n L: (h_n, j_n) \in H_z \times J_z, h_n
K\cap j_n L \neq \emptyset \}.
\end{equation}

We remark that it is possible that one set $h K$ may have nontrivial
intersections with two more sets $j_1 L$ and $j_2 L$, but $(h,j_1)$
and $(h,j_2)$ are counted differently in the set $\mathcal{A}$. Note
that $H_z K\cap J_z L \subset \cup \mathcal{A}$.

Define the set $\mathcal{B}=\{j_n^{-1}h_n: j_n^{-1}h_n K\cap L \neq
\emptyset, \; (h_n, j_n) \in H_z \times J_z \}$. We now show that
$\mathcal{B}$ is finite. Suppose not. By the convergence property,
there exists an infinite subsequence $\{j_{n_i}^{-1}h_{n_i}\}$ of
$\mathcal B$ such that $\{j_{n_i}^{-1}h_{n_i}\}$ converges locally
compactly to $b$ on $M \setminus \{a\}$, for some $a, b \in M$. We
claim $a =b$. Otherwise, using Lemma 2.5 in Bowditch \cite{Bow2}, we
have that $j_{n_i}^{-1}h_{n_i}$ are loxodromic elements for all
sufficiently large $n_i$. But this contradicts to the fact that
$\{j_{n_i}^{-1}h_{n_i}\}$ lie in the maximal parabolic subgroup
$G_z$.

Moreover, we have that $a=b=z$, since $z$ is the fixed point of
elements $j_n^{-1}h_n$. Note that $K\subset M\setminus \{z\}$ and
$L\subset M\setminus \{z\}$ are disjoint compact subsets. Since
$j_{n_i}^{-1}h_{n_i} K\cap L \neq \emptyset$, the subsequence
$\{j_{n_i}^{-1}h_{n_i}\}$ is a finite set by the convergence
property. This is a contradiction. Hence $\mathcal{B}$ is a finite
set.

Let $\mathcal{B}$ be a finite set, say $\{j_1^{-1}h_1,...,
j_r^{-1}h_r\}$, for example. Without loss of generality, we first
consider the elements in $\{j_n^{-1}h_n\}$ of the form $j_n^{-1}h_n
= j_1^{-1}h_1$.  Then $j_n j_1^{-1} = h_n h_1^{-1} \in H_z \cap J_z
= P$. We write $j_n = p_n j_1$ and $h_n = p_n h_1$ for some $p_n \in
P$. So we have $h_n K\cap j_n L = p_n (h_1 K\cap j_1 L)$ for each
$j_n^{-1}h_n = j_1^{-1}h_1$.

We can do the rewriting process similarly for other elements in
$\{j_n^{-1}h_n\}$, and finally we obtain $H_z K\cap J_z L \subset
\cup \mathcal{A} \subset P C$, where $C$ is a compact set defined as
$\bigcup\limits_{i=1}^r (h_i K \cap j_i L)$. The claim is proved.

Recall that we have proved there exists a compact subset $C \subset
M$ such that the following holds
\begin{equation}
\Lambda(H \cap J)\setminus \{z\} \subset (\Lambda(H) \cap
\Lambda(J))\setminus \{z\} \subset P C. \label{parabolic}
\end{equation}
We now have two cases to consider for finishing the proof of
proposition,

\item[\textbf{$P$ is finite}.] Since the right-hand of (\ref{parabolic}) is a compact set, there exists an open neighborhood of $z$ disjoint with $\Lambda(H)\cap\Lambda(J)$. Thus, $z$ is an isolated point of
$\Lambda(H)\cap\Lambda(J)$ and does not lie in $\Lambda(H \cap J)$.

\item[\textbf{$P$ is infinite}.] $P$ acts cocompactly on $\Lambda(H \cap J)\setminus \{z\}$. Thus, $z$ is a bounded parabolic point
of $H\cap J$.
\end{proof}

Summarizing the above results, we can now conclude with the proof of
Theorem \ref{thm:intersection}. Recall that by Theorem
\ref{thm:DQCRQC}, dynamically quasiconvex subgroups coincide with
relatively quasiconvex groups in relatively hyperbolic groups.
\begin{proof}[Proof of Theorem \ref{thm:intersection}]
By (QC-1) definition of relative quasiconvexity in \cite{Hru}, a
relatively quasiconvex subgroup acts on its limit set as a
geometrically finite convergence action. Then the limit set of a
relatively quasiconvex subgroup consists of conical limit points and
bounded parabolic points. Therefore, the decomposition of
$\Lambda(H) \cap \Lambda(J)$ follows from Propositions
\ref{prop:conical} and \ref{prop:boundpara}.
\end{proof}

\begin{rem}
In word hyperbolic groups, the exceptional set $E$ is empty since
there are no parabolic subgroups. In this case, limit sets of two
relatively quasiconvex subgroups intersect at least in two points,
once they intersect. But in the relative case, it is possible that
their limit sets intersect in only one (necessarily parabolic)
point. For example, let $H=\langle z+1\rangle, J = \langle
z+i\rangle$ be two parabolic subgroups of a Fuchsian group $G$
acting on the upper half space $\{z \in \mathbb{C}: \Im (z) > 0\}$.
Note that the intersection $H \cap J$ is trivial, but $H$ and $J$
share the same fixed point $\infty$.
\end{rem}

The following corollary follows from the isolatedness of the
exceptional set $E$.
\begin{cor}
Let $H, J$ be two relatively quasiconvex subgroups of $G$. If
$\Lambda(H) \subset \Lambda(J)$. Then either $\Lambda(H \cap J) =
\Lambda(H)$ or $H$ is a parabolic subgroup.
\end{cor}
\begin{proof}
Suppose $H$ is not a parabolic subgroup. Then $|\Lambda(H)| \ge 2$.
By Theorem \ref{thm:intersection}, we have $\Lambda(H) = \Lambda(H
\cap J) \sqcup E$, where $E$ consists of isolated points in
$\Lambda(H)$. It is well-known that limit sets are perfect, if
containing at least 3 points. So if $|\Lambda(H)| > 2$, then $E$ is
empty.

It suffices to consider the case $|\Lambda(H)| = 2$. In this case,
$H$ is a virtually cyclic group. Thus, $\Lambda(H)$ consists of two
conical limit points. By Lemma \ref{lem:conipara}, we have that $E$
is empty.
\end{proof}

\begin{center} \section[Proper Conjugates]{Proper Conjugates of Dynamically Quasiconvex Subgroups} \end{center}

Suppose $G$ has a convergence group action on $M$. According to
\cite{GMRS}, a subgroup $H \subset G$ is said to be \textit{maximal}
in its limit set if $H = Stab_G(\Lambda(H))$. Recall that Lemma
\ref{lem:samelimitset} shows any nonparabolic dynamically
quasiconvex subgroup is of finite index in the stabilizer of its
limit set.

\begin{lem} \label{lem:maxsubgrp}
Let $H$ be dynamically quasiconvex in a convergence group $G$ and
suppose $H$ is maximal in its limit set. Then for any $g \in G
\setminus H$, $g\Lambda(H) \nsubseteq \Lambda(H)$.
\end{lem}
\begin{proof}
If $|\Lambda(H)| = 1$, then $H$ is the maximal parabolic subgroup.
In this case, the conclusion is trivial. We now consider the case
$|\Lambda(H)| \geq 2$. By way of contradiction, we suppose that
$g\Lambda(H) \subseteq \Lambda(H)$.

Take distinct points $x$ and $y$ from $\Lambda(H)$. Since
$g\Lambda(H) \subseteq \Lambda(H)$, we have $g^n(x) \in \Lambda(H)$
and $g^n(y) \in \Lambda(H)$ for each $n \in \mathbb{N}$. Therefore,
we have that the cosets $g^{-n}H$ belong to the following set
\begin{center}
$\{gH \in G/H : g \Lambda(H) \cap \{x\} \neq \emptyset$ and $g
\Lambda(H) \cap \{y\} \neq \emptyset\}$.
\end{center}

Since $H$ is dynamically quasiconvex, we have the set $\{g^{-n}H \}$
is finite. Consequently there exist two different integers $m$ and
$n$ such that $g^{-m}H = g^{-n}H$, and thus $g^{n-m} \in H$. Then
$\Lambda(H) = g^{n-m} \Lambda(H) \subseteq g\Lambda(H) \subseteq
\Lambda(H)$. Hence $\Lambda(H) = g\Lambda(H)$, which is impossible
since $H$ is maximal in its limit set $\Lambda(H)$.
\end{proof}
\begin{rem}
Lemma \ref{lem:maxsubgrp} generalizes Lemma 2.10 in
Gitik-Mitra-Rips-Sageev \cite{GMRS}.
\end{rem}

We now prove Theorem \ref{thm:conjuagte}.

\begin{proof}[Proof of Theorem \ref{thm:conjuagte}]
Let $g$ be an element of $G \setminus H$ such that $gHg^{-1}
\subseteq H$. By Lemma \ref{lem:maxsubgrp}, it follows that $g$
belongs to the setwise stabilizer $K$ in $G$ of $\Lambda(H)$. Using
Lemma \ref{lem:samelimitset}, we obtain that $H$ is of finite index
in $K$. So $g^n$ belongs to $H$ for some $n$. Thus, we have $H = g^n
H g^{-n} \subset gHg^{-1} \subset H$. The proof is complete.
\end{proof}
\begin{rem}
The condition that $|\Lambda(H)| \ge 2$ could not be dropped. It is
known that there exists a finitely generated group $G$ containing a
finitely generated subgroup $H$ such that, for some $g \in G$,
$gHg^{-1} \subset H$ but $gHg^{-1} \subsetneq H$. See \cite{WZ} for
an elementary example. We then form a free product $G * F_2$, where
$F_2$ is a free group of rank 2. By the second definition of
relative hyperbolicity in \cite{Bow1}, $G * F_2$ is relatively
hyperbolic. In particular, $G$ is a maximal parabolic subgroup. But
$H \subset G * F_2$ doesnot satisfy the statement of Theorem
\ref{thm:conjuagte} for some $g \in G$.
\end{rem}

\begin{center} \section[Undistorted Subgroups]{Undistorted Subgroups of Groups with Nontrivial Floyd Boundary} \end{center}

In this section, we consider a finitely generated group $G$ with a
fixed finite generating set $S$, without assuming relative
hyperbolicity of $G$.

As usual, $S$ is assumed to be symmetric, i.e. $S = S^{-1}$. Then
the \textit{Cayley graph} $\GS$ of $G$ with respect to $S$, is
defined as an oriented graph with vertex set $G$ and edge set $G
\times S$. An edge $(g,s) \in G \times S$ goes from $g$ to $gs$.
Note that $\GS$ is a connected graph, which induces a \textit{word
metric} $d_S$ on $G$ by setting the length of each edge to be 1.

Given a rectifiable path $p$ in $\GS$, we denote by $p_-, p_+$ the
initial and terminal endpoint of $p$ respectively. Let $l(p)$ be the
length of $p$. We say $p$ is a \textit{$\epsilon$-quasigeodesic} for
a constant $\epsilon \ge 0$ if, for any subpath $q$ of $p$, we have
$l(q) < \epsilon d_S(q_-,q_+)+\epsilon$. Let $d_S(1, p)$ be the
distance from the identity to the path $p$ with respect to $d_S$.

Recall that a ($\epsilon$-)\textit{quasi-isometric} map $\phi : X
\to Y$ between two metric spaces $(X,d_X)$ and $(Y,d_Y)$ is a map
such that the following holds
\begin{center}
$\epsilon ^{-1} d_X (x, y) - \epsilon \leq d_Y (\phi(x), \phi(y))
\leq \epsilon d_X (x, y) + \epsilon$.
\end{center}

\begin{defn}
Let $H \subset G$ be a finitely generated subgroup with a finite
generating set $T$. Then $H$ is \textit{undistorted} if the
inclusion of $(H,d_T)$ into $(G,d_S)$ is a quasi-isometric map.
\end{defn}

Note that the definition of an undistorted subgroup is independent
of choices of finite generating sets $S, T$. Without loss of
generality, we assume that $T \subset S$ in the sequel. Then the
embedding $\imath: \HT \hookrightarrow \GS$ is a quasi-isometric
map. In particular, a geodesic in $\HT$ is naturally embedded as a
quasigeodesic in $\GS$.

We now briefly discuss the construction of Floyd boundary of a
finitely generated group. We refer the reader to \cite{Floyd},
\cite{Ka} and \cite{GePo1} for more details.

Let $f: \mathbb  N \to \mathbb{R}$ be the function $f(n)=n^{-2}$. We
rescale the length of each edge $e$ of $\GS$ by a factor
$f(d_S(1,e))$, and then take the Cauchy metric completion
$\overline{G}$. Denote by $\rho$ the complete metric on
$\overline{G}$. Then \textit{Floyd boundary} $\partial(G)$ is
defined as $\overline{G} \setminus G$. With a change of finite
generating sets, Floyd boundary is well-defined up to a bi-Lipschitz
homeomorphism.

If $\partial(G)$ consists of 0, 1 or 2 points then it is said to be
\textit{trivial}. Otherwise, it is uncountable and is called
\textit{nontrivial}. The class of groups with nontrivial Floyd
boundary includes non-elementary relatively hyperbolic groups
\cite{Ge2}.

In \cite{Ka}, Karlsson showed that if Floyd boundary is nontrivial,
then $G$ acts on $\partial(G)$ as a
convergence group action. In what follows, when speaking of limit
sets and dynamical quasiconvexity of subgroups in $G$, we have in
mind the convergence action of $G$ on $\partial(G)$.

The following lemma shows that the Floyd length of a far
(quasi)geodesic in $\GS$ is small. The original version was stated
in \cite{Ka} for geodesics, but its proof also works for
quasigeodesic in the Cayley graph.
\begin{lem} \label{lem:karlsson} \cite{Ka}
Given $\epsilon > 0$, there is a function $\Theta_\epsilon$:
$\mathbb N \to \mathbb R_{\ge 0}$ such that $\Theta_\epsilon(n) \to
0$ as $n \to \infty$ and the following property holds. Let $z, w$ be
two points in $G$ and let $\gamma$ be an $\epsilon$-quasigeodesic
between $z$ and $w$ of $\GS$. Then the following holds
\begin{center}
$\rho(z, w)  \leq \Theta_\epsilon(d_S(1,\gamma))$
\end{center}
\end{lem}

Recall that we assume $T \subset S$ and the embedding $\imath: \HT
\hookrightarrow \GS$ is quasi-isometric. The following lemma roughly
says that any two limit points of an undistorted subgroup $H$ can be
connected by a geodesic in $\HT$.

\begin{lem} \label{lem:geodesics}
If $H$ is undistorted in $G$ such that $\Lambda(H) \ge 2$, then
there exists a constant $\epsilon_0 \ge 0$ such that the following
holds. For any two distinct points $p, q \in \Lambda(H)$, there
exists an $\epsilon_0$-quasigeodesic $\gamma$ in $\GS$ between $p$
and $q$ such that $\gamma \subset \HT$.
\end{lem}

\begin{proof}
Since $p$ and $q$ are distinct limit points of $H$, there exist two
sequences $\{h_n\}$ and $\{h_n'\}$ of $H$ such that $h_n \rightarrow
p$ and $h_n' \rightarrow q$. Let $\delta = \rho(p,q)/3$.

Without loss of generality, we assume for all $n$, $h_n \in
B_\delta(p)$ and $h_n' \in B_\delta(q)$, after passage to
subsequences of $\{h_n\}$ and $\{h_n'\}$ respectively. Here,
$B_\delta(p)$ and $B_\delta(q)$ denote open metric balls centered at
$p$ and $q$ in $\overline{G}$ with radius $\delta$ respectively. It
then follows by the triangle inequality that
$\rho(h_n,h_n')>d(p,q)/3$ for all $n$.

Taking geodesics $\gamma_n$ in the Cayley graph $\HT$ such that
$(\gamma_n)_- = h_n$ and $(\gamma_n)_+ =h_n'$. By the
undistortedness of $H$, there is a positive constant $\epsilon_0$
depending on $H$, such that any geodesic in $\HT$ is an
$\epsilon_0$-quasigeodesic in $\GS$. Thus, $\gamma_n$ are
$\epsilon_0$-quasigeodesics in $\GS$. Observe that the endpoints
$h_n, h_n'$ of $\gamma_n$ have at least a $\rho$-distance $\delta$
in $\overline{G}$.

Let $\Theta_{\epsilon_0}$ be the function given by Lemma
\ref{lem:karlsson}. Since $\Theta_{\epsilon_0}(n) \to 0$ as $n \to
\infty$, let $R$ be the maximal integer $m$ such that
$\Theta_{\epsilon_0}(m) \ge \delta$. By Lemma \ref{lem:karlsson},
each quasigeodesic $\gamma_n$ intersects a closed ball $B$ centered
at identity with radius $R$ in $\GS$.

Therefore, using a Cantor diagonal argument based on $\gamma_n$, we
obtain an $\epsilon_0$-quasigeodesic $\gamma$ in $\GS$ between $p$
and $q$ such that the vertex set of $\gamma$ lies in $H$.
\end{proof}

\begin{rem}
In contrast with hyperbolic groups, two (quasi)geodesics in $\GS$
with same endpoints may not be uniformly Hausdorff distance bounded.
Thus, we could not guarantee that any (quasi)geodesic between $p$ and
$q$ satisfies the statement of Lemma \ref{lem:geodesics}.
\end{rem}

We are now ready to prove Theorem \ref{thm:dynconvex}.
\begin{proof}[Proof of Theorem \ref{thm:dynconvex}]
It suffices to establish the conclusion under the assumption
$|\Lambda(H)|\geq 2$. We are going to bound the following set
\begin{center}
$\{gH \in G/H : g \Lambda(H) \cap L \neq \emptyset$ and $g
\Lambda(H) \cap K \neq \emptyset\}$,
\end{center}
whenever $K$ and $L$ are disjoint closed subsets of $\partial(G)$.

Suppose, to the contrary, there exists a sequence of distinct cosets
$g_n H$ such that $g_n\Lambda(H) \cap K \neq \emptyset$ and
$g_n\Lambda(H) \cap L \neq \emptyset$. Let $p_n \in g_n\Lambda(H)
\cap K$ and $q_n \in g_n\Lambda(H) \cap L$. Note that
$g_n^{-1}(p_n), g_n^{-1}(q_n) \in \Lambda(H)$. By Lemma
\ref{lem:geodesics}, we obtain $\epsilon_0$-quasigeodesics
$\gamma_n$ between $g_n^{-1}(p_n)$ and $g_n^{-1}(q_n)$, such that
the vertex set of $\gamma_n$ lies in $H$. Hence $g_n(\gamma_n)$ are
$\epsilon_0$-quasigeodesics with endpoints $p_n, q_n \in g_n
\Lambda(H)$, such that the vertex set of $g_n (\gamma_n)$ lies in
the same coset $g_n H$.

Note that $\{p_n , q_n\} \in K \times L$. Since $K \times L$ is
compact in $\partial(G) \times \partial(G)$. there exists a uniform
positive constant $\mu$ depending on $K$ and $L$, such that
$\rho(p_n,q_n) \geq \mu$ for all $n$. Let $\Theta_{\epsilon_0}$ be
the function given by Lemma \ref{lem:karlsson}. Since
$\Theta_{\epsilon_0}(n) \to 0$ as $n \to \infty$, let $R$ be the
maximal integer $m$ such that $\Theta_{\epsilon_0}(m) \ge \mu$.

By Lemma \ref{lem:karlsson}, any $\epsilon_0$-quasigeodesic between
$p_n$ and $q_n$ intersects non-trivially with $B$, where $B$ is the
closed ball at the identity with radius $R$ in $\GS$. Let $c_n$ be an
intersection point of $g_n \gamma_n \cap B$. Then we have
$d_S(1,c_n) < R$ for every $n$.

Recall that the vertex set of $g_n(\gamma_n)$ lies in $g_n H$.
Therefore, for each $n$, there exists $h_n \in H$ such that $d_S(g_n
h_n, c_n)<1$. Then $d_S(1, g_nh_n)<R+1$ for all $n$. Since $S$ is a
finite set, we have the set $\{g_n h_n\}$ is finite. This is a
contradiction, as $\{g_n H\}$ is assumed as a sequence of different
$H$-cosets in $G$. The proof is complete.
\end{proof}

In view of Theorem \ref{thm:dynconvex}, the previous Lemma
\ref{lem:samelimitset}, Corollaries 2.7, 2.8 and Theorem
\ref{thm:conjuagte} can be stated in the setting of finitely
generated groups with nontrivial Floyd boundary. In favor of
applications in group theory, we state the following.
\begin{cor}
Let $H$ be undistorted in $G$ such that $|\Lambda(H)| \geq 2$. Then
$H$ is of finite index in its commensurator.  In particular, $H$ is
of finite index in its normalizer.
\end{cor}

\begin{cor}
Let $H$ be undistorted in $G$ such that $|\Lambda(H)| \geq 2$. Then
for any $g \in G \setminus H$, $gHg^{-1} \subseteq H$ implies that
$gHg^{-1} = H$.
\end{cor}

In relatively hyperbolic groups, the limit set of a relatively
quasiconvex subgroup consists of conical points and bounded
parabolic points. This fact allows us to complete the limit set
intersection theorem \ref{thm:intersection} for relatively
quasiconvex groups. In general convergence groups it is an
interesting question to ask whether dynamically quasiconvex
subgroups act geometrically finitely on their limit sets. The
following example gives a negative answer to the question.

\begin{example}
Dunwoody's inaccessible group $J$ in \cite{Dun} has infinite ends
and thus nontrivial Floyd boundary. Since infinitely ended groups
are relatively hyperbolic, by using a theorem of Stalling with the
second definition of relative hyperbolicity in \cite{Bow1}. Thus,
$J$ is relatively hyperbolic. Let $\mathbb H$ be a set of
representatives of conjugacy class of maximal parabolic subgroups of
$J$. Since $(J, \mathbb H)$ is relatively hyperbolic, then each $H
\in \mathbb H$ is undistorted in $J$. See the remark
\ref{rem:undistorted}.

In \cite{YANG3}, we prove that if $(G, \mathbb P)$ is relatively
hyperbolic, then $G$ acts geometrically finitely on $\partial G$ if
and only if each $P \in \mathbb P$ acts geometrically finitely on
its limit set $\Lambda(P) \subset \partial G$.

However, it is noted in \cite{YANG3} that $J$ doesnot act
geometrically finitely on its Floyd boundary $\partial J$.
Therefore, there exists at least one $H \in \mathbb H$ such that $H$
doesnot act geometrically finitely on its limit set $\Lambda(H)
\subset
\partial J$. Moreover, $\Lambda(H)$ contains infinite limit points,
otherwise $J$ would act geometrically finitely on $\partial J$.

In a word, $H$ is dynamically quasiconvex by Theorem
\ref{thm:dynconvex}, but doesnot act geometrically finitely on its
limit set in the Floyd boundary $\partial J$.
\end{example}

\bibliographystyle{amsplain}

\end{document}